\documentclass[11pt]{amsart}

\newtheorem{teo}{Theorem}[section]
\newtheorem{prop}[teo]{Proposition}
\newtheorem{cor}[teo]{Corollary}
\newtheorem{lema}[teo]{Lemma}
\newtheorem{claim}[teo]{Claim}

\theoremstyle{definition}
\newtheorem{dfn}[teo]{Definition}

\theoremstyle{remark}
\newtheorem{rem}[teo]{Remark}

\numberwithin{equation}{section}

\newcommand{\dist}{\ensuremath{\mathrm{dist} }}

\newcommand{\tub}{\ensuremath{\mathrm{Tub} }}

\newcommand{\F}{\ensuremath{\mathcal{F}}}
\newcommand{\TO}{\ensuremath{\mathcal{O}}}
\newcommand{\singularF}{\ensuremath{\mathcal{X}_{F}}}

\newcommand{\dank}{\textsf{Acknowledgments:\ }}

\begin{document}

\title{Proofs of conjectures about singular riemannian foliations }


\author{Marcos M. Alexandrino} 
 
\address{Marcos Martins Alexandrino\\
Instituto de Matem\'{a}tica e Estat\'{\i}stica,
Universidade de S\~{a}o Paulo (USP),
Rua do Mat\~{a}o, 1010, Bloco A,
05508 090, S\~{a}o Paulo, Brazil}

\email{malex@ime.usp.br}

\email{marcosmalex@yahoo.de}


\date{September, 2005.}


\begin{abstract}

We prove that if the normal distribution of a singular riemannian foliation is integrable, then each leaf of this normal distribution  can be extended to be a complete immersed totally geodesic submanifold (called section), which meets every leaf orthogonally. In addition the set  of regular points is open and dense in each section. 
This result generalizes a result of Boualem and solves a problem inspired by a remark of Palais and Terng and a work of Szenthe about polar actions.

We also study the singular holonomy of a singular riemannian foliation with sections (s.r.f.s. for short) and in particular the tranverse orbit of the closure of each leaf. Furthermore we prove that the closure of the leaves of a s.r.f.s on $M$ form a partition of $M$ which is a singular riemannian foliation. This result proves partially a conjecture of Molino.   

\end{abstract}

\maketitle

\section{Introduction}

In this section we will recall some definitions 
and  state our main results as Theorem \ref{teo-PT-conjecture}
and Theorem \ref{teo-Molino-conjecture}. 

We start by recalling the definition of a singular riemannian foliation (see the book of P. Molino \cite{Molino}).

\begin{dfn}
 A partition $\F$ of a complete riemannian manifold $M$ by connected immersed submanifolds (the \emph{leaves}) is called a \emph{singular foliation} of $M$ if it verifies condition (1) and {\it singular riemannian foliation} if it verifies condition (1) and (2):

\begin{enumerate}
\item $\F$ is \emph{singular},
i.e., the module $\singularF$ of smooth vector fields on $M$ that are tangent at each point to the corresponding leaf acts transitively on each leaf. In other words, for each leaf $L$ and each $v\in TL$ with foot point $p,$ there is $X\in \singularF$ with $X(p)=v$.
\item  The partition is \emph{transnormal}, i.e., every geodesic that is perpendicular at one point to a leaf remains perpendicular to every leaf it meets.
\end{enumerate}
\end{dfn}

Let $\F$ be a singular riemannian foliation on a complete riemannian manifold $M.$  A leaf $L$ of $\F$ (and each point in $L$) is called \emph{regular} if the dimension of $L$ is maximal, otherwise $L$ is called {\it singular}. Let $L$ be an immersed submanifold of a riemannian manifold $M.$  A section $\xi$ of  the normal bundle $\nu L$ is said  to be a \emph{parallel normal field} along $L$ if  $\nabla^{\nu}\xi\equiv 0,$ where $\nabla^{\nu}$ is the normal connection.  $L$  is said to have a {\it globally flat normal bundle},  if the holonomy  of the normal bundle  $\nu L$ is trivial, i.e., if any normal vector can be extended to a globally defined parallel normal field. 

Now we recall the concept of singular riemannian foliation with section, which was introduced in 
 \cite{Alex1}, \cite{Alex2}  and continued to be studied by  T\"{o}ben and I in \cite{Toeben2}, \cite{Alex3}, \cite{AlexToeben} and recently by  Lytchak and  Thorbergsson in \cite{LytchakThorbergsson}.  

\begin{dfn}[s.r.f.s.]
 Let $\F$ be a singular riemannian foliation on a complete riemannian manifold $M.$
$\F$ is said to be a \emph{singular riemannian foliation with sections} (s.r.f.s. for short) if for each regular point $p,$ the set $\Sigma :=\exp_{p}(\nu_p L_{p})$ is a complete immersed submanifold that meets each leaf orthogonally. $\Sigma$ is called a \emph{section}.
\end{dfn}
\begin{rem}
In \cite{Alex2} it was also assumed that the set of regular points is open and dense in $\Sigma.$ As we will see in Theorem \ref{teo-PT-conjecture} this is a consequence of the definition above. We will explain in Remark \ref{rem-s.r.f.s-equifocal}
 the importance of the density of the set of regular points in a section.
\end{rem}

Let $p\in M$ and $\tub(P_{p})$ be a tubular neighborhood of a plaque $P_{p}$ that contains $p.$ Then the connected component of $\exp_{p}(\nu P_{p})\cap \tub(P_{p})$  that contains $p$ is called a \emph{slice} at $p.$ Let $S_{p}$ denote it. 
Now consider a section $\Sigma,$ a point $p\in\Sigma$  and a convex normal neighborhood $B$ of $p$ in $M.$ Then a connected component of $B\cap\Sigma$ that contains $p$ is called \emph{local section} (centered at $p)$.   
These two concepts  are  related to each other. Indeed,  the slice at a singular point is the union of the local sections that contain this singular point (see Theorem \ref{sliceteorema}).

A typical example of a s.r.f.s is the partition formed by the orbits of a polar action. 

An isometric action of a compact Lie group $G$ on a complete riemannian manifold $M$ is called \emph{polar} if there exists a complete immersed submanifold $\Sigma$ of $M$ 
that meets all $G$-orbits orthogonally and whose dimension is equal to the codimension of the regular $G$-orbits. We call $\Sigma$ a \emph{section}. A typical example of a polar action is a compact Lie group with a biinvariant metric that acts on itself by conjugation. In this case the maximal tori are the sections. One can find a large number of examples of polar action on  symmetric spaces in the work of Kollross \cite{Ko}and Podest\`{a} and Thorbergsson \cite{PoTh}.

In \cite{PTlivro}[page 96, Remark 5.6.8] Palais and Terng proposed a conjecture that can be formulated as follows: \emph{ Let $G$ be an isometric action of a compact Lie group on $M$ such that the distribution of the normal space to the regular orbit is integrable. Then there exists a complete totally geodesic immersed section for the action of $G$ which meets all orbits perpendicularly.} In \cite{HOL} Heintze, Olmos and Liu proved that Palais and Terng were right in their conjecture. In particular they proved that the set of regular points is dense in each section (see  Tebege \cite{Samuel} for an alternative proof of the density of regular points in a section).
We would like to stress that Szenthe \cite{Szenthe} worked on the above problem even before Terng and Palais. In fact, Szenthe proved, under some extra assumptions, the existence of a generalized section, i.e., a submanifold (not necessarily complete) that intersect all the leaves orthogonally. 


In \cite{Boualem} Boualem proved  the  following result:
   \emph{Let $\F$ be a singular riemannian foliation such that the closure of each leaf is compact and the distribution of normal spaces to the regular leaves  is integrable. Then 
each leaf of this normal distribution  can be extended to be an immersed totally geodesic submanifold (called generalized section), which meets every leaf orthogonally.} 

He did not prove that the so called generalized section is complete and hence could not give a complete answer to the conjecture of Palais and Terng for the case of singular riemannian foliation. Also the question of the density of the set of regular points in a section was not considered.

In this paper we prove the following result, which gives a complete answer to  the conjecture of Palais and Terng for the case of  singular riemannian foliations.

\begin{teo}
\label{teo-PT-conjecture}
Let $\F$ be a singular riemannian foliation on a complete manifold $M.$ Suppose that the distribution of normal spaces of the regular leaves  is integrable. Then $\F$ is s.r.f.s.
In addition  the set of regular points is open and dense in each section. 
\end{teo}

We note that Lytchak  and Thorbergsson \cite{LytchakThorbergsson} have recently  proved this result under the additional condition that singular riemannian foliation has no horizontal conjugate point.

Before we state the other main result of this paper, we must briefly recall the notion of isoparametric and equifocal submanifolds. These concepts are fundamental to understand local and global aspects of the theory of s.r.f.s. 
 
A submanifold of an euclidean space is called \emph{isoparametric} if its normal bundle is flat and if the principal curvatures along any parallel normal vector field are constant. The  singular partition formed by the parallel submanifolds of an isoparametric submanifold in an euclidean space  is a typical example of s.r.f.s. The history of isoparametric submanifolds and their generalizations can be found in the survey \cite{Th} of  Thorbergsson (see also \cite{Th2}). 

According to the slice theorem (see Theorem \ref{sliceteorema}) the restriction of a s.r.f.s to a slice at a singular point  is a s.r.f.s which is diffeomorphic to an isoparametric foliation on an open set of an euclidean space. Hence the slice theorem give us a description of the plaques of a s.r.f.s. on a riemannian manifold. However, it does not guarantee that different plaques belong to different leaves. To obtain this kind of information, we need the concept of singular holonomy which is related to the concept of local equifocal submanifold reviewed below.

 In \cite{TTh1}  Terng and  Thorbergsson introduced the concept of equifocal submanifold with  flat sections in symmetric space in order to generalize the definition of isoparametric submanifold in Euclidean space. Now we review the definition of equifocal submanifolds in a riemannian manifold. 

\begin{dfn} 
\label{dfn-equifocal}
A connected immersed submanifold  $L$  of a complete Riemannian manifold $M$ is called \emph{equifocal} if
\begin{enumerate}
\item[0)] the normal bundle $\nu(L)$ is globally flat,
\item[1)] $L$ has sections, i.e., for each $p\in L$ 
 the set $\Sigma :=\exp_{p}(\nu_p L_{p})$ is a complete, immersed, totally geodesic submanifold.
\item[2)] For each  parallel normal field $\xi$  along $L,$  the derivative of  the map    $\eta_{\xi}:L\rightarrow M,$ defined as $\eta_{\xi}(x):=\exp_{x}(\xi),$ has constant rank.
\end{enumerate}
\end{dfn}
A connected immersed submanifold $L$ is called \emph{locally equifocal} if, for each $p\in L,$ there exists a neighborhood $U\subset L$ of $p$ in $L$ such that $U$ is an equifocal submanifold.

\begin{rem}[Relation between s.r.f.s. and equifocal submanifolds]
\label{rem-s.r.f.s-equifocal}

In \cite{TTh1} Terng and  Thorbergsson proved that the singular partition formed by the parallel submanifolds of an equifocal submanifold $L$ on $M$ is a s.r.f.s, when $M$ is a symmetric space of compact type and the sections are flat. In \cite{Alex2} I proved that the leaves of a s.r.f.s. on a riemannian manifold $M$ are locally equifocal (see  Theorem \ref{frss-eh-equifocal}). In this proof the property that the set of regular points is dense in a section plays an important role.
In \cite{Toeben2} T\"{o}ben  used the blow up technique to study local equifocal submanifolds (which he called submanifold with parallel focal structure). 
He gave a necessary and sufficient condition for a closed embedded local equifocal submanifold to induce a s.r.f.s. 
 T\"{o}ben's Theorem (for the case of equifocal submanifold) follows from the next result, which is proved in \cite{Alex3}. \emph{
Let $L$ be a closed embedded equifocal submanifold.  Suppose that set of regular points of each section $\Sigma$ is an open and dense set in $\Sigma.$ Define  $\Xi$ as the set of all parallel normal fields along $L.$ Then $\F:=\{\eta_{\xi}(L)\}_{\xi\in \, \Xi}$ is a s.r.f.s.
} Here  a point of a section is called regular if there exists only one local section  that contains $p,$ i.e., if given two local sections $\sigma$ and $\tilde{\sigma}$  that contain $p,$ they have the same germ at $p.$

\end{rem}

The property that each leaf of a s.r.f.s is locally equifocal allow us to extend the holonomy map to the whole local section and hence to define  the notion of singular holonomy map (see Proposition \ref{prop-holonomia-singular}).  If we  fix a local section $\sigma,$ we can define a Weyl pseudogroup of $\sigma$ as the pseudogroup generated by the  singular holonomy maps defined on $\sigma$ (see Definition \ref{prop-holonomia-singular}). The Weyl pseudogroup of a local section centered at $q$  and the s.r.f.s on a slice $S_{q}$ describe completely the s.r.f.s. on a neighborhood of $q.$     

We are finally ready to recall the conjecture that motivate our next result.

In \cite{Molino} Molino proved that, if $M$ is compact, the closure of the leaves of a (regular) riemannian foliation form a partition of $M$ which is a singular riemannian foliation. He also proved that the leaf closure are orbits of a locally constant sheaf of germs of (transversal) Killing fields.  

If the foliation is singular riemannian foliation and $M$ is compact, then Molino was able to prove (see \cite{Molino} Theorem 6.2 page 214) that the closure of the leaves should be a transnormal system, but as he remarked, it remains to prove that the closure of the leaves is in fact a singular foliation. In this paper we prove the Molino's conjecture, when $\F$ is a s.r.f.s. In addition we study the singular holonomy of $\F$ and in particular the   tranverse orbits of the closure of a leaf. We will not assume that $M$ is compact.

\begin{teo}
\label{teo-Molino-conjecture}
Let $\F$ be a s.r.f.s. on a complete riemannian manifold $M.$ 
Then
\begin{enumerate}
\item[a)] the closure of the leaves of $\F$ form a partition of $M$ which is a singular riemannian foliation,i.e, $\{\overline{L}\}_{L\in\F}$ is a singular riemannian foliation.
\item[b)] Each point $q$ is contained in an homogenous submanifold $\TO_{q}$ (possible with dimension $0).$ If we fix a local section $\sigma$ that contains $q,$ then $\TO_{q}$ is a connected component of an orbit of the closure of the Weyl pseudogroup of $\sigma.$   
\item[c)] If $q$ is a point of the submanifold $\overline{L},$ then a neighborhood of $q$ in $\overline{L}$ is the product of the homogenous submanifold  $\TO_{q}$ with plaques with the same dimension of the plaque $P_{q}.$
\item[d)] Let $q$ be a singular point and  $T$  the intersection of the slice $S_{q}$ with the   singular stratum that contains $q.$ Then the normal connection of $T$ in $S_{q}$ is flat.
\item[e)] Let $q$ be a singular point and $T$ defined as in Item d). Let $v$ be a parallel normal vector field along $T,$ $x\in T$ and $y=\exp_{x}(v)$. Then  
$\TO_{y}=\eta_{v}(\TO_{x}).$ 
\end{enumerate}
\end{teo}

One can construct examples that illustrate the above theorem by means of suspension of homomorphisms (see \cite{Alex2} for details). In fact, the suspension technique is very useful to construct examples of s.r.f.s. with nonembedded leaves, with exceptional leaves and also inhomogeneous examples. Other techniques to construct examples of s.r.f.s on nonsymmetric spaces are  suitable changes of metric and  surgery (see \cite{AlexToeben} for details). 

\begin{cor}
\label{cor-estrato-trivial-folhas-fechadas}
Let $\F$ be a s.r.f.s. on a complete manifold $M$ and $q$ a singular point. Let $T$ denote the intersection of the slice $S_{q}$ with the  stratum that contains $q.$ Suppose that $T=\{q\}.$ Then all the leaves of $\F$ are closed.
\end{cor}

\begin{rem}
According to the slice theorem (see Theorem \ref{sliceteorema})  the restriction of the foliation $\F$ to the slice $S_{q}$ is diffeomorphic to an isoparametric foliation $\widetilde{\F}$ on an open set of an euclidean space. Therefore the condition that $T$ is a point is equivalent to saying that a regular leaf of $\widetilde{\F}$ is a full isoparametric submanifold.
\end{rem}

\begin{rem}
\label{rem-quase-polar}

By the suspension technique one can construct an example of s.r.f.s.  such that, for a local section $\sigma,$ the isometric action of 
$\overline{W}_{\sigma}$ (the closure of the Weyl pseudogroup of $\sigma$) is not polar. In particular,  this implies  that there exists a s.r.f.s $\F$ such that the partition formed by the closure of the leaves of $\F$ is not a s.r.f.s. 
Although the isometric action of  $\overline{W}_{\sigma}$  need not to be polar in a neighborhood of a singular point $q,$ if $\sigma$ is flat,  there always exists  a totally geodesic submanifold $N_{q}$ that is orthogonal to the orbits of $\overline{W}_{\sigma}.$
 This submanifold $N_{q}$ is a slice of $T$ in $\sigma$ (see Equation \ref{def-slice-N}).
It fails to be a section of $\overline{W}_{\sigma}$, if its  dimension is lower than the codimension of the regular orbits of $\overline{W}_{\sigma}.$ To prove that $N_{q}$ is orthogonal to the orbits of $\overline{W}_{\sigma}$ which it meets, one can use 
the above theorem and a  result of Heintze, Olmos and Liu\cite{HOL}.

\end{rem}

This paper is organized as follows. In Section 2 we review some facts about s.r.f.s. and fix the notation. In Section 3 and 4  we  prove Theorem \ref{teo-PT-conjecture} and Theorem \ref{teo-Molino-conjecture} respectively. In Section 5 we prove Corollary \ref{cor-estrato-trivial-folhas-fechadas}. 

\dank I am grateful to  Professor Gudlaugur Thorbergsson and Professor Claudio Gorodski for  useful suggestions.

\section{Facts about s.r.f.s.}

In this section we recall some results about s.r.f.s. that  will be used in this work. Details can be found in \cite{Alex2}. Throughout this section we assume that $\F$ is a s.r.f.s. on a complete riemannian manifold $M.$

Let us start with a result that relates s.r.f.s. to equifocal submanifolds (see Definition \ref{dfn-equifocal} to recall the definitions of equifocal submanifolds and the endpoint map $\eta_{\xi}$).

\begin{teo}
\label{frss-eh-equifocal}
The regular leaves of $\F$ are  equifocal. In particular, the union of regular leaves that 
have trivial normal holonomy is an open and dense set in $M,$ provided that all the leaves are compact.
\end{teo}

A consequence of the first statement of the theorem is that given a regular leaf $L$ with trivial holonomy then we can reconstruct $\F$ by taking all parallel submanifolds of $L.$ More precisely we have

\begin{cor}
\label{cor-map-paralelo}
Let $L$ be a regular leaf of  $\F.$ 
\begin{enumerate}
\item[a)] Let  $\beta$ be a smooth curve of $L$ and  $\xi$ a parallel normal field along  $\beta$. Then the curve $\eta_{\xi}\circ \beta$ belongs to a leaf of \F.
\item[b)] Let  $L$ be a regular leaf with trivial holonomy and $\Xi$ denote the set of all parallel normal fields along $L.$ Then $\F=\{\eta_{\xi}(L)\}_{\xi\in \, \Xi}.$ 
In particular, if $\xi$ is a parallel normal field along $L$ then the endpoint map $\eta_{\xi}:L\rightarrow L_{q}$ is surjective, where $q=\eta_{\xi}(x)$ for $x\in L.$
\end{enumerate}
\end{cor}

Corollary \ref{cor-map-paralelo} allows us to define singular holonomy map, which will be very useful to study $\F.$  

\begin{prop}[Singular Holonomy]
\label{prop-holonomia-singular}
Let $L_{p}$ be a regular leaf,  $\beta$ a smooth curve in $L_{p}$ and let $[\beta]$ denote the homotopy class of $\beta.$ Let $U$ be a local section centered at $p=\beta(0).$ Then there exists a local section $V$ centered at $\beta(1)$ and 
an isometry $\varphi_{[\beta]}:U\rightarrow V$  that has the following properties:
\begin{enumerate}
\item[1)]$\varphi_{[\beta]}(x)\in L_{x}$ for each $x\in U,$
\item[2)]$d\varphi_{[\beta]}\xi(0)=\xi(1),$ where $\xi$ is a  parallel normal field along $\beta.$
\end{enumerate}
\end{prop}
An isometry as above is called \emph{singular holonomy map along $\beta$}. 
 
We remark that, in the definition of the singular holonomy map, singular points can be contained in the domain $U.$  
If the domain $U$ and the range $V$ are sufficiently small, then the singular holonomy map coincides with the holonomy map along $\beta.$

Now we  recall some results about the local structure of $\F$, in particular about the structure of the set of singular points in a local section.


Note that, if $q$ is a singular point, the restriction $\F|S_q$ of $\F$ to the slice $S_{q}$ is also a singular foliation. In fact, since $\singularF$ acts transitively on the leaves,  the plaques (the connected components of the leaves intersected with $\tub(P)$) are transversal to $S_q.$ 

The relation between a slice $S_{q},$ a local section and $\F|S_q$ are given by the next result. 

\begin{teo}[Slice Theorem]
\label{sliceteorema}
Let $q$ be a singular point of $M$ and $S_{q}$ a slice at $q.$ Then
\begin{enumerate}
\item[a)] Let $\epsilon$ be the radius of the slice $S_{q}.$ Denote  $\Lambda(q)$  the set of 
local sections $\sigma$ containing $q,$ such that $\dist(p,q)<\epsilon$ for each $p\in\sigma.$ 
 Then  $S_{q}= \cup_{\sigma\in\Lambda (q)}\, \sigma.$
\item[b)] $S_{x}\subset S_{q}$ for all $x\in S_{q}.$
\item[c)] $\F|S_q$ is a s.r.f.s. on $S_{q}$ with the induced metric of $M.$
\item[d)] $\F|S_q$ is diffeomorphic  to an isoparametric foliation on an open set of $\mathbf{R}^{n},$  where $n$ is the dimension of   $S_{q}.$
\end{enumerate}
\end{teo}

 From d) it is not difficult to derive the following corollary.

\begin{cor}
\label{estratificacao-singular}
Let $\sigma$ be a local section. Then the set of singular points of $\F$ contained in $\sigma$ is a finite union of totally geodesic hypersurfaces. 
These hypersurfaces are sent by a diffeomorphism to focal hyperplanes 
 contained in a section of an isoparametric foliation on an open set of a euclidean space.   
\end{cor}
We will call the set of singular points of $\F$ contained in $\sigma$ the \emph{singular stratum of the local section} $\sigma$. 
Let $M_{r}$ denote the set of regular points in $M.$ A \emph{Weyl Chamber} of a local section $\sigma$ is the closure in $\sigma$ of a connected component of $M_{r}\cap\sigma.$  

\begin{cor}
\label{convexidadedeWeyl}
A Weyl Chamber of a local section is a convex set.
\end{cor}

\begin{cor}
\label{pontos-singulares-na-geodesica}
Let $\gamma$ be a geodesic orthogonal to a regular leaf. Then the  singular points are isolated on $\gamma.$ 
\end{cor}

\begin{rem}
The above fact can be proved without the Slice Theorem. In fact it was used in \cite{Alex2}  to prove the Slice Theorem.
\end{rem}

The Slice Theorem establishes a relation between s.r.f.s. and isoparametric foliations. By analogy with the classical theory of isoparametric submanifolds, it is 
 natural to ask if we can define a (generalized) Weyl group action on $\sigma.$ The next definitions and results answer this question.

\begin{dfn}[Weyl Pseudogroup $W$]
\label{definitionWeylPseudogroup}
 The pseudosubgroup  generated by all singular holonomy  maps $\varphi_{[\beta]}$ such that $\beta(0)$ and $\beta(1)$ belong to the same local section $\sigma$ is called \emph{generalized Weyl pseudogroup of} $\sigma.$ Let $W_{\sigma}$ denote this pseudogroup. In a similar way we define $W_{\Sigma}$ for a section $\Sigma.$ Given a  slice $S$ we will define $W_{S}$ as the set of  all singular holonomy maps   $\varphi_{[\beta]}$ such that $\beta$ is contained in the slice $S.$
\end{dfn}

\begin{rem}
To recall the definition of pseudogroups and  orbifolds  see E. Salem \cite[Appendix D]{Molino}. 
\end{rem}

\begin{prop}
\label{propWisinvariant}
Let $\sigma$ be a local section. Then the reflections in the hypersurfaces of the singular stratum of the local section $\sigma$ leave $\F|\sigma$  invariant. Moreover these reflections are elements of $W_{\sigma}.$
\end{prop}

One can construct an example of a s.r.f.s. by suspension such that  $W_{\sigma}$ is larger than the pseudogroup generated by the reflections  in the  hypersurfaces of the singular stratum of $\sigma.$ The next result gives a sufficient condition 
to guarantee that both pseudogroups coincide.

\begin{prop}
\label{Weylgrouparereflection}
Suppose that each leaf of $\F$ is compact and has trivial normal holonomy. Let $\sigma$ be a local section.  Then  $W_{\sigma}$ is generated by the reflections  in the hypersurfaces of the singular stratum of the local section. 
\end{prop}

\begin{rem}
In \cite{AlexToeben} T\"{o}ben and I proved that if $M$ is simply connected and the leaves of $\F$ are compact, then each regular leaf has trivial holonomy.
\end{rem}

Finally let us recall a  result related to Theorem \ref{teo-Molino-conjecture}.

\begin{prop}
\label{estrutura-transversa-do-fecho} 
Let $\sigma$ be a local section. Consider a point  $p\in\sigma,$ then
\[\overline{W_{\sigma}}\cdot p= \overline{L_{p}\cap\sigma}.\]
In other words, the closure of $L_{p}\cap\sigma$ is an orbit of 
complete closed pseudogroup of local isometries. In particular $\overline{L_{p}\cap\sigma}$ is a closed submanifold.
\end{prop}

\section{Proof of  Theorem \ref{teo-PT-conjecture}}

In this section we will prove Theorem \ref{teo-PT-conjecture}.

A singular riemannian foliation is called \emph{integrable} if  the distribution of the normal spaces of the regular leaves  is integrable. A neighborhood of a leaf of this  distribution is called \emph{regular section}. 

In what follows we will need the following result due to Boualem (see Lemma 1.2.1, Proposition 1.2.3 and Lemma 1.2.4 of \cite{Boualem}). 

\begin{prop}[Boualem \cite{Boualem}]
\label{prop-Boualem}
Let $\F$ be an integrable singular riemannian foliation on a complete riemannian manifold $M.$ Let $S_{q}$ denote the slice at $q$ of radius $\epsilon.$
Let $\sigma$ be a regular section contained in  the normal ball $B_{\epsilon}(q).$ Suppose also that $\sigma$  intersects the slice $S_{q}.$ 
Then
\begin{enumerate}
\item[a)]
 There exists a subspace $V\subset T_{q}S_{q}$ such that $\sigma\subset D_{\epsilon}(q):=\exp_{q}(V\cap U),$ 
where $U$ is the ball of radius $\epsilon$  in $T_{q}S_{q}$ centered at $0.$
  \item[b)] the map $\exp^{-1}_{q}$ is a diffeomorphism between   the foliation $\F|_{S_{q}}$ and an integrable singular riemannian foliation   $\widetilde{\F}$ on the open set $U$ of the inner product space  $(T_{q}S_{q},<,\,>_{q}),$ where  $<,\,>_{q}$ denote the metric of $T_{q}M.$ 
Furthermore, the submanifold  $\widetilde{\sigma}:=\exp_{q}^{-1}(\sigma)$  is a regular section of $\widetilde{\F}.$
\end{enumerate} 
\end{prop}

\begin{rem}
\label{rem-prop-Boualem}
Let $<,>_{0}$ denote the canonical euclidean product. Then there exists a positive  definite symmetric matrix $A$ such that $<X,Y>_{q}=<A\,X,Y>_{0}$  The isometry $\sqrt{A}:(T_{q} S_{q},<,>_{q})\rightarrow (\mathbf{R}^{n},<,>_{0})$ is a diffeomorphism between  the foliation $\widetilde{\F}$ and an integrable singular riemannian foliation  $\widehat{\F}$ on an open set of the inner product space $(\mathbf{R}^{n},<,>_{0}).$
\end{rem}

\begin{prop}
\label{pontos-singulares-na-geodesica}
Let $\F$ be an integrable singular riemannian foliation on a complete riemannian manifold $M$ and $\gamma$ be a geodesic orthogonal to a regular leaf. Then the set of singular points is isolated on $\gamma.$ 
\end{prop}
\begin{proof} Since the set of regular points on $\gamma$ is open, we can suppose that $q=\gamma(0)$ is a singular point and that $\gamma(t)$ is a regular point for $-\delta<t<0.$ We will show that there exists $\epsilon>0$ so that the point  $\gamma(t)$ is  a regular point if  $0<t<\epsilon.$

It follows from Proposition \ref{prop-Boualem} that there exist an open set $U\subset T_{q}S_{q}$ and  an integrable singular riemannian foliation $\widetilde{F}$ on $U$  such that  $\exp^{-1}_{q}(\F)=\widetilde{\F}.$ Set $\widetilde{\gamma}:=\exp^{-1}_{q}(\gamma).$

Note that, for $t_{0}<0,$ the point  $0=\widetilde{\gamma}(0)$ is a focal point of the regular leaf $\widetilde{L}_{\widetilde{\gamma}(t_{0})}.$ Indeed, since $0$ is a leaf of $\widetilde{\F}$ with  dimension zero and $\widetilde{\F}$ is a singular riemannian foliation, each geodesic that joins $0$ to  $\widetilde{L}_{\widetilde{\gamma}(t_{0})}$ must be orthogonal to $\widetilde{L}_{\widetilde{\gamma}(t_{0})}.$ 

 Since focal points are isolated along $\widetilde{\gamma},$ we can choose $\epsilon>0,$ such that $ \widetilde{\gamma}(t)$ is not a focal point of the plaque  $\widetilde{P}_{\widetilde{\gamma}(t_{0})}$ along $\widetilde{\gamma}$ for $0<t<\epsilon.$ 

Suppose that there exists $0<t_{1}<\epsilon$ such that $x=\widetilde{\gamma}(t_{1})$ is a singular point. 

Let $\sigma$ be a regular section of $\widetilde{\gamma}(t_{0}).$ It follows from Proposition \ref{prop-Boualem} that $\sigma$ is contained in a subspace $V$ of dimension $k.$ Note that $x\in V.$ 
Now we choose straight lines $l_{1},\ldots,l_{k}:[0,1]\rightarrow V$ such that 
\begin{enumerate}
\item $l_{i}(0)\in\sigma,$
\item $l_{i}(1)=x,$
\item $\{l_{i}^{'}(1)\}$ is a basis to $T_{x}V.$
\end{enumerate}

Since $\sigma$ is a regular section and $\widetilde{\F}$ is a singular riemannian foliation, then the straight lines $l_{i}$ are orthogonal to the plaque $\widetilde{P}_{x}$. Therefore we have: 
\[\sigma \subset\nu_{x}:=x+\nu_{x}\widetilde{P}_{x}.\]
Since $x$ is a singular point, it follows that:
\[\dim \sigma<\dim \nu_{x}.\]

The equations above imply that $\dim  \widetilde{P}_{\widetilde{\gamma}(t_{0})}\cap\nu_{x}>0.$ Thus we can find geodesics in the subspace $\nu_{x}$ that join $x$ to the submanifold  $ \widetilde{P}_{\widetilde{\gamma}(t_{0})}\cap\nu_{x}.$ Since the foliation is  riemannian, these geodesics are also orthogonal to $ \widetilde{P}_{\widetilde{\gamma}(t_{0})}\cap\nu_{x}.$ This implies that  $x$ is a focal point of the submanifold $\widetilde{P}_{\widetilde{\gamma}(t_{0})}$. This contradicts our choice of $\epsilon$ and completes the proof. \end{proof}


\begin{lema}
\label{prop-F-isoparametrica}
Let $S_{q}$ be the slice at $q\in M.$ 
Let  $\widetilde{\F}$ be the integrable singular riemannian foliation on an open set $U$ of $T_{q}S_{q}$  defined in Proposition \ref{prop-Boualem}. Consider a subspace $V$ orthogonal to a regular leaf of $\widetilde{\F}.$ 
Then
\begin{enumerate}
\item[a)]  The set of regular points of $\widetilde{\F}$ is open and dense in $V\cap U,$ which is  orthogonal to the leaves of $\widetilde{\F}.$
In particular the foliation $\widehat{\F}$ defined in Remark \ref{rem-prop-Boualem} is an isoparametric foliation.
\item[b)]the set of regular points of $\F$ is open and dense in $D_{\epsilon}(q):=\exp_{q}(V\cap U).$
\item[c)]The submanifold $D_{\epsilon}(q)$ is orthogonal to the leaves of $\F.$
\item[d)] $D_{\epsilon}(q)$ is a totally geodesic submanifold.
 \end{enumerate}
\end{lema}
\begin{proof}
a) By assumption there is a regular section $\widetilde{\sigma}\subset V.$  The fact that the set of regular points of $\widetilde{\F}$ is open and dense in $V\cap U$ follows direct from Proposition \ref{pontos-singulares-na-geodesica}. To prove that $V\cap U$ is orthogonal to the plaque $\widetilde{P}_{x}$ at $x\in V\cap U,$ define
straight lines $l_{1},\ldots,l_{k}:[0,1]\rightarrow V$ such that 
\begin{enumerate}
\item $l_{i}(0)\in\widetilde{\sigma},$
\item $l_{i}(1)=x,$
\item $\{l_{i}^{'}(1)\}$ is a basis to $T_{x}V.$
\end{enumerate}
Since $\widetilde{\sigma}$ is a regular section and $\widetilde{\F}$ is a singular riemannian foliation,  the straight lines $l_{i}$ are orthogonal to the plaque $\widetilde{P}_{x}$. To conclude the proof of item a) we note that the properties of $V$ imply that $\widetilde{\F}$ is a s.r.f.s. Now the slice theorem (Theorem \ref{sliceteorema}) implies that $\widehat{\F}$ is an isoparametric foliation.

b) The item b) follows direct from item a).

c) Let $y\in D_{\epsilon}(q)$ be a regular point and $\sigma_{y}$ a regular section that contains $y.$ We want to prove that $\sigma_{y}\subset D_{\epsilon}(q).$ Define $\widetilde{\sigma}_{y}:=\exp_{q}^{-1}(\sigma_{y}).$ It follows from Proposition \ref{prop-Boualem} that  $\widetilde{\sigma}_{y}$ is a regular section of $\widetilde{\F}.$ Now item a) implies that $\widetilde{\sigma}_{y}\subset V\cap U.$ Therefore  $\sigma_{y}\subset D_{\epsilon}(q),$ i.e., the regular leaf $L_{y}$ is orthogonal to $D_{\epsilon}(q)$ at $y.$  Since regular points are dense in $D_{\epsilon}(q),$ the singular leaves are also orthogonal to $D_{\epsilon}(q).$ 

d) It follows from item c) and from the fact that regular points are dense in $D_{\epsilon}(q).$ 
\end{proof}

\begin{prop}
\label{cor-discos-geodesicos}
Let $q$ be a point of $M$ and $D_{\epsilon}(q)$ the submanifold  defined in Proposition \ref{prop-Boualem}. Let $B_{\delta}(q)$ be the normal ball of radius $\delta$ in $M$ centered at $q.$ Define the disc $D_{\delta}(q):=\{\exp_{q}(v)$ such that $v\in T_{q}D_{\epsilon}(q)$ and $\|v\|<\delta\}.$  Then
\begin{enumerate}
\item[a)]the set of regular points of $\F$ is open and dense in $D_{\delta}(q).$
\item[b)]The submanifold $D_{\delta}(q)$ is orthogonal to the leaves of $\F.$
\item[c)] $D_{\delta}(q)$ is a totally geodesic submanifold.
\end{enumerate}
\end{prop}
\begin{rem}
Note that $\delta$ (the radius of the normal neighborhood of $q)$   is often greater then  $\epsilon$ (the radius of the  slice $S_{q}$) which depends on the leaf $L_{q}.$ 
\end{rem}
\begin{proof}
Let $\xi\in T_{q}D_{\epsilon}(q)$ be a vector such that the geodesic $\gamma(t):=\exp_{q}(t\xi)$ has regular points ($0\leq t\leq 1$). It follows from Proposition \ref{pontos-singulares-na-geodesica} that the singular points are isolated on $\gamma.$ Since $\gamma$ is compact, we can cover this geodesic with a finite number of slices with centers $\gamma(t_{i})$ so that $\gamma(t_{i})$ is the unique possible singular point in $\gamma\cap S_{\gamma(t_{i})}.$ Now we can construct a submanifold $T$ by gluing the totally geodesic submanifolds $D_{\epsilon}(\gamma(t_{i}))$ defined in Lemma \ref{prop-F-isoparametrica}. Note that $T$ satisfies the conditions of item a) b) and c). 

Since $\exp_{q}$ is a diffeomorphism and $T$ is totally geodesic, we can find a neighborhood $W$ of $\xi$ in $T_{q} D_{\epsilon}(q)$ and define a neighborhood of $\gamma$ in $T$ as $N:=\exp_{q}(W).$ Therefore $N\subset D_{\delta}(q).$ Hence we have constructed a neighborhood $N$ of $\gamma$ in $D_{\delta}(q)$ that satisfies the conditions of item a), b) and c).

To conclude the proof, note that we can choose the vector $\xi$ arbitrarily close to any vector of $T_{q}D_{\epsilon}(q)$ for the set of regular points of $\F$  is an open and dense set in $D_{\epsilon}(q).$ 
  
\end{proof}

\textbf{Proof of the Theorem}

Let $q\in M$ be a point.
 Define $D_{r_{1}}(q):=\{\exp_{q}(v)$ $ / v\in T_{q} D_{\epsilon}(q)$ and  $\|v\|\leq r_{1}\},$ where $r_{1}$ is a positive number. 
Since $ D_{r_{1}}(q)$ is compact, we can find a compact neighborhood of $ D_{r_{1}}(q)$ and a number $\delta>0$ such that, for each point $x$ in this neighborhood, $B_{\delta}(x)$ is a normal ball. Let $\gamma$ be a geodesic with length $r_{1}$ such that  $\gamma^{'} (0)\in T_{p} D_{\epsilon}(q).$ 

At first we note that we can extend $D_{\epsilon}(q)$ along the geodesic $\gamma$ to a noncomplete submanifold which satisfies the items of Proposition \ref{cor-discos-geodesicos}
, gluing totally geodesic discs of radius $\delta/2.$ To see this let the first disc be $D_{\delta/2}(q).$ It follows from Proposition \ref{cor-discos-geodesicos} that $D_{\delta/2}(q)$ satisfies  the  items a),b) and c) . Now we define the second disc as 
$D_{\delta/2}(\gamma(\delta/4))$ $:=\{\exp_{\gamma(\delta/4)}(v)$ for $v\in T_{\gamma(\delta/4)}D_{\delta/2}(q)$ and $ \|v\|\leq \delta/2\}.$ Again, it follows from Proposition \ref{cor-discos-geodesicos} that  $D_{\delta/2}(\gamma(\delta/4))$ satisfies  the items a), b) and c). We continue in this way, defining discs $D_{i}$'s by induction and extending $D_{\epsilon}(q)$ along $\gamma.$   

Let $\Sigma_{r_{1}}$ denote the union of the extensions constructed above along every geodesics $\gamma$ with length $r_{1}$ such that $\gamma (0)^{'}\in T_{q}D_{\epsilon}(q).$

\begin{lema}
$\Sigma_{r_{1}}$ can be covered by a finite number of discs $D_{i}$'s (of radius $\delta$) which satisfy the items a), b) and c) of Proposition \ref{cor-discos-geodesicos}.
\end{lema}

\begin{proof}

We start the proof with the following claim.
   
\begin{claim}
\label{claim-intersecao-finita}
Let $x\in D_{r_{1}}(q).$ Then  $D_{r_{1}}(q)\cap B_{\delta/4}(x)$  belongs to a finite union of discs $D_{i}$'s (of radius $\delta/2$) which satisfy the items a), b) and c) of Proposition  \ref{cor-discos-geodesicos}.    
\end{claim}


To check this claim,  let $K$ denote the set $\{v\in T_{p}D_{\epsilon}(q)$ such that $ \exp_{p}(v)\in \overline{B_{\delta/4}(x)}$ and $ \|v\|\leq r_{1}\}$ and choose a vector  $v_{0}\in K.$ 
As explained above,  we can extend $D_{\epsilon}(q)$ along the geodesic $\exp_{p}(t\,v_{0})$ to  a totally geodesic noncomplete submanifold  $T.$  Since $\exp$ is continuous, we can find a neighborhood  $W$ of the segment   $v_{0},$ such that $\exp_{q}(t\,u)\subset T$  for all   $u\in W\cap K$   and    $ 0<t\leq 1.$   In other words,  we can find a neighborhood $W$ of $v_{0}$ such that for all $u\in W\cap K,$  the extensions along  $\exp_{q}(t\,u)$ are the same in  $\overline{B_{\delta/4}(x)}.$ The claim now follows  from the fact that $K$ is compact.  

To conclude the proof of the lemma, we cover $D_{r_{1}}(q)$ with a finite number of balls $B^{i}_{\delta/4}.$ Then, using Claim \ref{claim-intersecao-finita}, we deduce that 
$D_ {r_{1}}(q)$ is contained in a finite  union of discs $D_{\delta/2}(p_{i})$ which satisfy
the items a), b) and c) of Proposition \ref{cor-discos-geodesicos}. Thus, by the definition of $\Sigma_{r_{1}},$ we conclude that,  for each $x\in \Sigma_{r_{1}}$, there exists $i_{0}$ such that $x\in D_{\delta}(p_{i_{0}}).$    


\end{proof}
 
We are ready to show  that $\Sigma_{r_{1}}$ is an immersed submanifold in $M.$ Let $\mathbf{G}^{k}(M)$ be the $k$-Grassmann bundle of $T M$ and $\Pi:\mathbf{G}^{k}(M)\rightarrow M$  the natural  projection.   Let  $D^{\star}_{i}\subset \mathbf{G}^{k}(M)$  be the set $\{T_{x}D_{i}/$ $ x\in D_{i}\}$ and  $\Sigma^{\star}_{r_{1}}\subset \mathbf{G}^{k}(M)$  be the set $\{T_{x}D_{i}/$ for all $x\in\tilde{\Sigma}_{r_{1}}\}.$  Hence we have defined  the blow-up's of    $D^{\star}_{i}$ and $\Sigma_{r_{1}}.$ It is easy to see that $D^{\star}_{i}$ is embedded in $\mathbf{G}^{q}(M).$  In addition, if  $D^{\star}_{i}$ and  $D^{\star}_{j}$ have a common point then they coincide in a neighborhood of this point.  Since $\Sigma^{\star}_{r_{1}}$  can be covered by a finite numbers of $D^{\star}_{i},$  we see that     $\Sigma^{\star}_{r_{1}}$  is  embedded  in $\mathbf{G}^{q}(M).$  We conclude that $\Sigma_{r_{1}}$ is immersed for $\Sigma_{r_{1}}=\Pi\Sigma^{\star}_{r_{1}}.$ 

To finish the proof, we choose a sequence $r_{1}< \cdots < r_{n}\rightarrow \infty $  and construct $\Sigma_{r_{i}}$ and $\Sigma^{\star}_{r_{i}}$ in a similar way as we have done before.     
As above we can also show that $\Sigma^{\star}_{r_{i}}$ is embedded and $\Sigma_{r_{i}}$ is immersed. Finally we define $\Sigma^{\star}:=\cup_{r_{i}}\Sigma^{\star}_{r_{i}}.$  We can see that  $\Sigma^{\star}$ is a submanifold   and that $\Sigma:=\cup_{r_{i}}\Sigma_{r_{i}}$  is immersed for $\Sigma=\Pi \Sigma^{\star}.$   
The Hopf and Rinow's Theorem implies that $\Sigma$ is complete. 




\section{Proof of Theorem \ref{teo-Molino-conjecture}}

In this section we will prove some propositions and claims, which  will imply direct the Theorem \ref{teo-Molino-conjecture}. In fact, item a) of the theorem  will follow from Corollary \ref{cor-particao-transnormal} and Proposition \ref{prop-fechoF-eh-singular}. Item b) will follow from Proposition \ref{prop-orbitas-W} and Definition \ref{dfn-tranverse-orbits}. The items c) and d)  will be proved in Proposition \ref{prop-produto-OT-plaqueta} and Proposition \ref{conexao-singular-stratum-flat} respectively. Finally, item e)  will follow from Claim \ref{prop-fechoF-eh-singular-claim3}.

\begin{prop}
\label{prop-orbitas-W}
Let $S_{q}$ denote a slice at a point $q\in M$ and consider two local sections $\sigma_{0}$ and $\sigma_{1}$ that contain $q.$ Then
\begin{enumerate}
\item[a)] There exists $h\in W_{S_{q}}$ such that $W_{\sigma_{1}}=h W_{\sigma_{0}} h^{-1}.$
\item[b)] $\overline{W_{\sigma_{1}}}\cdot q =\overline{W_{\sigma_{0}}} \cdot q \subset \sigma_{1}\cap \sigma_{0}.$
\end{enumerate} 
\end{prop}
\begin{proof}

a)Let $w\in W_{\sigma_{1}}.$ Then there exists a curve $\beta_{1}$ contained in a regular leaf such that $\beta_{1}(0)\in\sigma_{1},$ $\beta_{1}(1)\in\sigma_{1}$ and $w=\varphi_{[\beta_{1}]}.$

Let $\alpha$ be a curve contained in a regular leaf of $\F\cap S_{q}$ such that $\alpha(0)\in\sigma_{0}$ and $\alpha(1)\in\sigma_{1}.$ Set $h:=\varphi_{[\alpha]}.$

Now we define two curve that  join the local section $\sigma_{0}$ to $\sigma_{1}$ as 
 $\rho(t):=\varphi_{[\alpha|_{[0,t]}]}(h^{-1}\beta_{1}(0))$ and   $\varrho(t):=\varphi_{[\alpha|_{[0,t]}]}(h^{-1}\beta_{1}(1)).$ Finally we define the concatenation
$\beta_{0}:=\varrho^{-1}*\beta_{1}*\rho.$ 
Since $h=\varphi_{[\alpha]}=\varphi_{[\rho]}=\varphi_{[\varrho]},$  we have that $\varphi_{[\beta_{0}]}=h^{-1}\varphi_{[\beta_{1}]}h.$  
From the last equation, we conclude that $W_{\sigma_{1}}\subset  h W_{\sigma_{0}} h^{-1}.$ By  a similar argument one can prove that 
$W_{\sigma_{1}}\supset h W_{\sigma_{0}} h^{-1}$ and conclude the proof of item a).

b)If $q$ is a regular point, then $\sigma_{0}=\sigma_{1}$ and the result is trivial.
Hence we can assume that $q$ is a singular point. 

We can deduce the next claim from the Slice Theorem. 
\begin{claim}
\label{prop-orbitas-W-claim}
Let $T$ be the intersection of the local section $\sigma_{0}$ with the singular stratum that contains $q.$ Then $T\subset \sigma_{0}\cap\sigma_{1}$ and $h|_{T}$ is the identity.
\end{claim}
 On the other hand, it follows direct from the definition of $W_{\sigma_{0}}$ that 
\begin{eqnarray}
\label{prop-orbitas-W-eq-1}
W_{\sigma_{0}}\cdot q\subset T.
\end{eqnarray}

Now  Equation  \ref{prop-orbitas-W-eq-1} and Claim \ref{prop-orbitas-W-claim}
imply that
\begin{eqnarray}
\label{prop-orbitas-W-eq-2}
h W_{\sigma_{0}}\cdot q= W_{\sigma_{0}}\cdot q.
\end{eqnarray}
Finally,  using  Proposition \ref{estrutura-transversa-do-fecho}, 
Item a), Claim \ref{prop-orbitas-W-claim} and Equation \ref{prop-orbitas-W-eq-2}, we can prove Item b) as follows:

\begin{eqnarray*}
\overline{W_{\sigma_{1}}}\cdot q &\stackrel{\mathrm{Prop.}\ref{estrutura-transversa-do-fecho}} {=}& \overline{W_{\sigma_{1}}\cdot q} \\
   &\stackrel{\mathrm{Item\ a)}}{=}&\overline{h W_{\sigma_{0}}h^{-1}\cdot q}\\
    &\stackrel{\mathrm{Claim\ }\ref{prop-orbitas-W-claim}}{=}&\overline{h W_{\sigma_{0}}\cdot  q}\\
     &\stackrel{\mathrm{Eq. \ }\ref{prop-orbitas-W-eq-2}}{=}&\overline{W_{\sigma_{0}}\cdot q}\\                               
  &\stackrel{\mathrm{Prop.}\ref{estrutura-transversa-do-fecho}} {=}& \overline{W_{\sigma_{0}}} \cdot q 
\end{eqnarray*}

\end{proof}

\begin{dfn}[Transverse Orbits]
\label{dfn-tranverse-orbits}
Consider a point $q$ of $M$ and $\sigma$ a local section that contains $q.$ Let $\TO_{q}$ denote the submanifold $ \overline{W}_{\sigma}^{0} \cdot q,$ where $ \overline{W}_{\sigma}^{0}$ is the connected component of the identity element. It follows from Proposition \ref{prop-orbitas-W} that  this definition is independent of the choice of the local section $\sigma.$ 
\end{dfn}


\begin{prop}
\label{prop-produto-OT-plaqueta}
Let $L$ be a leaf of $\F$ and $\overline{L}$ the closure of $L.$ Then $\overline{L}$ is a submanifold of $M.$
In addition if $q$ is a point of $\overline{L},$ then a neighborhood of $q$ in $\overline{L}$ is the product of $\TO_{q}$ with plaques with the same dimension of $P_{q}.$
\end{prop}

\begin{proof}

Let $S_{q}$ be a slice at $q.$

\begin{claim}
\label{prop-produto-OT-plaqueta-claim1}
Let $T$ denote the intersection of $S_{q}$ with the  stratum that contains $L_{q}.$ Then $L\cap S_{q}\subset T.$  
\end{claim}
To check the above claim, we can assume that $q$ is a singular point, because the claim is trivial when $q$ is a regular point. Suppose that  $L\cap S_{q}$ is not contained in $T.$ Then Proposition \ref{prop-holonomia-singular}  implies that $L\cap S_{q}$ must be equidistant from $T.$ Thus  $q$ does not belong to $\overline{L}$ what  contradicts  our assumption.

From Claim \ref{prop-produto-OT-plaqueta-claim1} we can deduce the next claim.
\begin{claim}
\label{prop-produto-OT-plaqueta-claim2}
$L\cap S_{q}$ are points of $S_{q}.$  
\end{claim}
Now Claim \ref{prop-produto-OT-plaqueta-claim2} and the item b) of the slice theorem (Theorem \ref{sliceteorema}) imply the following claim.
\begin{claim}
\label{prop-produto-OT-plaqueta-claim3}
If $x\in\overline{L}\cap S_{q},$ then 
\begin{enumerate}
\item the plaque $P_{x}\subset\overline{L},$ 
\item the plaque $P_{x}$ intersect $S_{q}$ in one point,
\item the plaque $P_{x}$ is orthogonal to $S_{q}.$
\end{enumerate}
\end{claim}
According to Claim \ref{prop-produto-OT-plaqueta-claim3}, to prove the proposition it suffices  
to show that the connected component of $\overline{L}\cap S_{q}$ which contains $q$ is the submanifold $\TO_{q}.$ 

At first we note that if  $q$ is a regular point, then $S_{q}$ is a local section and it follows from Proposition \ref{estrutura-transversa-do-fecho} 
that the connected component of $\overline{L}\cap S_{q}$ is the submanifold $\TO_{q}.$

Now we suppose that $q$ is a singular point and choose a local section $\sigma$ that contains $T.$ 
Since $q\in\overline{L},$ we have that $q\in \overline{W}_{\sigma}\cdot z,$ for $z\in L.$
This implies that $\overline{W}_{\sigma}\cdot q=\overline{W}_{\sigma}\cdot z.$ Therefore $\TO_{q}$ is the connected component of $\overline{W}_{\sigma}\cdot z$ that  contains $q.$
To conclude the proof we note that  Proposition \ref{estrutura-transversa-do-fecho} implies that   the connected component of $\overline{W}_{\sigma}\cdot z$ that  contains $q$
is the connected component of $\overline{L}\cap S_{q}$ that contains $q.$

\end{proof}

\begin{cor}
\label{cor-produto-OT-plaqueta}
Let $L$ be a leaf of $\F,$ $\overline{L}$ the closure of $L$ and $q\in\overline{L}.$ Then $L_{q}\subset \overline{L}$ and $\overline{L}_{q}=\overline{L}.$
\end{cor}

\begin{cor}
\label{cor-particao-transnormal}
The singular partition $\{\overline{L}\}_{L\in\F}$ is a transnormal system.
\end{cor}
\begin{proof}
Let $\gamma$ be a geodesic orthogonal to $\overline{L}$ so that $\gamma(0)\in \overline{L}.$  Since $L_{\gamma(0)}\subset \overline{L}$
(see Corollary \ref{cor-produto-OT-plaqueta}), the geodesic $\gamma$ is orthogonal to the leaf $L_{\gamma(0)}.$ Therefore, by the slice theorem, the segment of geodesic $\gamma|_{(-\epsilon,\epsilon)}$ is containded in a local section $\sigma,$ for a small $\epsilon.$

Since $\gamma$ is orthogonal to $\overline{L},$ it follows from  Proposition \ref{prop-produto-OT-plaqueta} that $\gamma$ is orthogonal to $\overline{W}_{\sigma}\cdot\gamma(0).$
Since $\overline{W}_{\sigma}$ is a pseudogroup of isometries, $\gamma$ is orthogonal to the orbits of $\overline{W}_{\sigma}.$ This fact and Proposition \ref{prop-produto-OT-plaqueta} imply that $\gamma$ is orthogonal to $\overline{L}_{\gamma(t)}$ for every $t\in (-\epsilon,\epsilon).$   

With the above argument one can prove that  $\gamma|_{(-r,r)}$ is orthogonal to $\{\overline{L}\}_{L\in\F}.$ To prove that $\gamma$ is also orthogonal to $\overline{L}_{\gamma(r)},$ we  note that $\gamma$ is orthogonal to the leaf $L_{\gamma(r)}$ and hence is contained in a local section $\sigma.$ By Proposition \ref{prop-produto-OT-plaqueta},  $\gamma$ is orthogonal to the orbits of $\overline{W}_{\sigma}$ and in particular to the orbit $\overline{W}_{\sigma}\cdot\gamma(r).$ Thus $\gamma$ is orthogonal to $\overline{L}_{\gamma(r)}.$ The result follows then from the connectedness of $\gamma.$

\end{proof}

\begin{prop}
\label{conexao-singular-stratum-flat}
Let $q$ be a singular point and  $T$  the intersection of the slice $S_{q}$ with the  stratum that contains $q.$ Then the normal connection of $T$ in $S_{q}$ is flat. 
\end{prop}
\begin{proof}
Let $x_{0}$ and $x_{1}$ be points of $T$ and $v$ be a vector of the normal space $\nu_{x_{0}}T\cap T_{x_{0}}S_{q}.$
Let $||^{\nu}(\alpha)_{0}^{1}v$ denote the parallel normal  transport of $v$ along a curve $\alpha\subset T$ so that $\alpha(0)=x_{0}$ and $\alpha(1)=x_{1}.$ We want to prove that 
 $||^{\nu}(\alpha)_{0}^{1}v$ independ of the choise of the curve $\alpha.$
Let $||(\alpha)_{0} ^{t} v$ denote the parallel transport of $v$ along $\alpha.$ Since $T$ is totally geodesic, $||^{\nu}(\alpha)_{0} ^{1} v = ||(\alpha)_{0}^{1}v.$ Thus, it suffices to prove that $||(\alpha)_{0}^{1}v$ independ of the choise of the curve $\alpha$ that joins $x_{0}$ to $x_{1}.$

We choose a local section $\sigma$ that is tangent to $v.$ It follows from the slice theorem that $T\subset \sigma.$ Since $\sigma$ is totally geodesic, $|| (\alpha)_{0} ^{t} v\subset \sigma$ for all $t.$

From now on we consider $\sigma$ as the ambient space of $T.$ 

Let $N_{x_{0}}$ be the slice  at $x_{0}$ of $T$  i.e.,  
\begin{eqnarray}
\label{def-slice-N} 
N_{x_{0}}:=\{\exp_{x_{0}}v, \ \mathrm{for}\ v\in(\nu_{x_{0}}T\cap T_{x_{0}}\sigma) \ \mathrm{and} \ \exp_{x_{0}}v\subset S_{q}\}.
\end{eqnarray}

It follows from the slice theorem that we can choose a basis $\{\xi_{i}\}_{i=1\cdots k}$ of $T_{x_{0}}N_{x_{0}}$ such that $x_{0},\exp_{x_{0}}(\xi_{1}),\cdots,\exp_{x_{0}}(\xi_{k})$ are vertices of  a simplicial cone contained in the singular stratification of $N_{x_{0}}.$ Let $T_{i}$ denote the singular stratum that contains the curve $s\rightarrow \exp_{x_{0}}(s\xi_{i}),$ for small $s.$

Since $T_{i}$ is totally geodesic, the parallel transport $||(\alpha)_{0}^{t}\xi_{i}$ is always tangent to $T_{i}.$ On the other hand, $T$ is a totally geodesic hipersurface of $T_{i}.$ These two facts implies that $||(\alpha)_{0}^{1}\xi_{i}$ independ of the choice of $\alpha.$ Therefore 
$||(\alpha)_{0}^{1} v$ independ of the choice of $\alpha$ because $\{\xi_{i}\}$ is a basis of $N_{x_{0}}.$ 

\end{proof}

\begin{prop}
\label{prop-fechoF-eh-singular} 
The singular partition $\overline{\F}:=\{\overline{L}\}_{L\in\F}$ is a singular foliation. 
\end{prop}
\begin{proof}

Let $S_{q}$ denote the slice at a point $q.$
It suffices to prove that $\overline{\F}\cap S_{q}$ is a singular foliation.
Since $\F\cap S_{q}$ is already a singular foliation, Proposition \ref{prop-produto-OT-plaqueta}
implies that it is enough to prove that $\{\TO_{x}\}_{x\in S_{q}}$ is a singular foliation.

If $q$ is a regular point, then the slice $S_{q}$ is a local section $\sigma$  and the result follows from the fact that  $\{\TO_{x}\}_{x\in \sigma}$ are  orbits of the pseudogroup of isometries $\overline{W}_{\sigma}.$ Thus we can suppose that $q$ is a singular point.

Let  $T$ denote  the intersection of the slice $S_{q}$ with the  singular stratum that contains $q.$ Since the normal connection of $T$ in $S_{q}$ is flat (see Proposition \ref{conexao-singular-stratum-flat}), we can define the end point map $\eta_{\xi}:T\rightarrow S_{q}$ where $\xi$ is a parallel normal vector field along $T.$ 
One can easily check that:
\begin{claim}
\label{prop-fechoF-eh-singular-claim1} 
The partition $\{\eta_{\upsilon}(T)\}_{\upsilon\in\Upsilon}$ is a regular foliation, where $\Upsilon$ is the set of all parallel normal vector field along $T.$
\end{claim}

Using Claims \ref{prop-fechoF-eh-singular-claim1} above and Claim \ref{prop-fechoF-eh-singular-claim3} below 
one can show that  
$\{\TO_{x}\}_{x\in S_{q}}$ is a singular foliation and conclude the proof.

\begin{claim}
\label{prop-fechoF-eh-singular-claim3}
Let $x\in T$ and $y=\exp_{x}(v)$ for $v\in \nu_{x}T.$ Then 
$\TO_{y}=\eta_{v}(\TO_{x})$ 
\end{claim}

We start the proof of Claim \ref{prop-fechoF-eh-singular-claim3} choosing a 
 local section $\sigma$ so that $v$ is tangent to $\sigma.$ 
We also consider a curve $w:(-\epsilon,\epsilon)\rightarrow \overline{W}^{0}_{\sigma}$ in the connected component of the identity of $\overline{W}_{\sigma}$ and define a curve $\alpha(t):=w(t)\cdot x.$   
Next we review the definitions used  in the proof of Proposition \ref{conexao-singular-stratum-flat}. Let $N_{\alpha(t)}$ be the slice at $\alpha(t)$ of $T$ (see  (\ref{def-slice-N}))  
 and  $\{\xi_{i}\}$ a basis of $T_{x}N_{x}$ so that 
$x,\exp_{x}(\xi_{1}),\cdots,\exp_{x}(\xi_{k})$ are vertices of  a simplicial cone contained in the singular stratification of $N_{x}.$ Let $T_{i}$ denote  the singular stratum that contains the curve $s\rightarrow \exp_{x}(s\, \xi_{i}),$ for small $s.$

Since the local isometries of $\overline{W}_{\sigma}$ leave $T$ invariant, we can conclude that 
\begin{eqnarray}
\label{prop-fechoF-eh-singular-eq-0}
d\, w(t)_{x}\, \xi_{i} \in T_{\alpha(t)}N_{\alpha(t)}
\end{eqnarray}

We can also deduce that $w(t):T_{i}\rightarrow T_{i},$ using the fact that  $\overline{W}_{\sigma}$ leave singular stratum invariant. 
 Therefore we  have:
\begin{eqnarray}
\label{prop-fechoF-eh-singular-eq-1} 
 d w(t)_{x}\,\xi_{i}\in T_{\alpha(t)}T_{i}.
\end{eqnarray}
On the other hand, the fact that $T_{i}$ is totally geodesic and $T$ is a totally geodesic hypersurface of $T_{i}$ imply that:
\begin{eqnarray}
\label{prop-fechoF-eh-singular-eq-2} 
||(\alpha)_{0}^{t}\xi_{i}\in T_{\alpha(t)}T_{i}\cap T_{\alpha(t)} N_{\alpha(t)},
\end{eqnarray}
where $||(\alpha)_{0}^{t}\xi_{i}$ is the parallel transport of $\xi_{i}$ along the curve $\alpha$ 

Finally we can deduce from Equations \ref{prop-fechoF-eh-singular-eq-0} , \ref{prop-fechoF-eh-singular-eq-1} and \ref{prop-fechoF-eh-singular-eq-2}  that:
\begin{eqnarray}
\label{prop-fechoF-eh-singular-eq-3} 
d w(t)_{x}\,\xi_{i}=||(\alpha)_{0}^{t}\xi_{i}.
\end{eqnarray} 
Since $\{\xi_{i}\}$ is a basis of $N_{x},$ it follows from Equation \ref{prop-fechoF-eh-singular-eq-3}  that 
\begin{eqnarray}
\label{prop-fechoF-eh-singular-eq-4} 
d w(t)_{x}\, v=||(\alpha)_{0}^{t}v,
\end{eqnarray} 
which implies Claim \ref{prop-fechoF-eh-singular-claim3}.

\end{proof}

\section{ Proof of Corollary \ref{cor-estrato-trivial-folhas-fechadas}}

Let $\sigma$ be a local section that contains $q.$ The fact that $T=\{q\}$ implies that $W_{\sigma}$ has discrete orbits. Therefore the leaves in a tubular neighborhood of $L_{q}$ are closed. Now suppose that there exists a nonclosed leaf $L.$ Let $\gamma$ be a geodesic orthogonal to $\overline{L}$ such that $\gamma(0)\in\overline{L}$ and $\gamma(1)$ is a regular point of $\sigma.$ Since $L$ is not closed, the submanifold $\TO_{\gamma(0)}$ has dimension greater then zero. It follows from item c) of Theorem \ref{teo-Molino-conjecture}
 that $\gamma$ is orthogonal to $\TO_{\gamma(0)}.$ Since  the leaves near to a neighborhood of $\gamma(1)$ in $\sigma$ intersect $\sigma$ only one time, we can use Corollary \ref{cor-map-paralelo} to conclude that all points of $\gamma$ near to $\gamma(1)$ are focal points of $\TO_{\gamma(0)}.$ This contradicts the fact that focal points are isolated.


\bibliographystyle{amsplain}

\end{document}